\documentclass[12pt]{article}
\usepackage{mathrsfs}
\usepackage{amsfonts}
\usepackage{stmaryrd}
\usepackage{latexsym,amssymb,amsmath}
\begin{document}

\baselineskip=18pt
\newcommand{\headingstobeshown}{}
\def\tenrm{\rm}

\renewcommand{\theequation}{\arabic{section}.\arabic{equation}}

%\makeatletter      % '@' is now a normal "letter" for TeX
%\@addtoreset{equation}{section}
%\makeatother       % '@' is restored as a "non-letter" character for TeX

\newenvironment{proof}
               {\begin{sloppypar} \noindent{\bf Proof.}}
               {\hspace*{\fill} $\square$ \end{sloppypar}}

\newtheorem{atheorem}{\bf \temp}[section]

\newenvironment{theorem}[1]{\def \temp{#1}
                \begin{atheorem}\medskip
                }
                {\medskip
                 \end{atheorem}}

\newtheorem{thm}[atheorem]{Theorem}
\newtheorem{cor}[atheorem]{Corollary}
\newtheorem{lem}[atheorem]{Lemma}
\newtheorem{pro}[atheorem]{Property}
\newtheorem{prop}[atheorem]{Proposition}
\newtheorem{de}[atheorem]{Definition}
\newtheorem{rem}[atheorem]{Remark}
\newtheorem{fac}[atheorem]{Fact}
\newtheorem{ex}[atheorem]{Example}
\newtheorem{pr}[atheorem]{Problem}
\newtheorem{cla}[atheorem]{Assert}

\begin{center}{\Large \bf  Cohomology of Oriented Tree Diagram Lie Algebras}\end{center}

\vspace{0.2cm}

\begin{center}{\large Li Luo }
\end{center}

\begin{center}{Institute of Mathematics, Academy of Mathematics \&
System Sciences,
%}\end{center}
%\begin{center}{
\\Chinese Academy of Sciences, Beijing 100190, China}\end{center}
\begin{center}{E-mail: luoli@amss.ac.cn}\end{center}

\vspace{0.4cm}
\begin{abstract}
Xu introduced a family of root-tree-diagram nilpotent Lie algebras
of differential operators, in connection with  evolution partial
differential equations. We generalized his notion to more general
oriented  tree diagrams. These algebras are natural analogues of the
maximal nilpotent Lie subalgebras of finite-dimensional simple Lie
algebras. In this paper, we use Hodge Laplacian to study the
cohomology of these Lie algebras. The ``total rank conjecture" and
``$b_2$-conjecture" for the algebras are proved. Moreover, we find
the generating functions of the Betti numbers by means of Young
tableaux for the Lie algebras associated with certain tree diagrams
of single branch point. By these functions and
Euler-Poincar$\acute{\mbox{e}}$ principle, we obtain analogues of
the denominator identity for finite-dimensional simple Lie algebras.
The result is a natural generalization of the Bott's classical
result in the case of  special linear Lie algebras.
\end{abstract}

\vspace{0.5cm}
\section{Introduction}

Cohomology of Lie algebras are important objects in mathematics,
which are related to the geometry of the corresponding Lie groups,
invariant differential operators, combinatorial identities,
integrable systems,  Riemannian foliations and cobordism theory [F].
In particular, cohomology of nilpotent Lie algebras with
coefficients in the trivial module are more commonly used and has
attracted many mathematicians' attention. However, there are only a
few results on the full cohomology of Lie algebras up to this stage.

Santharoubane [S] found the cohomology of Heisenberg Lie algebras.
Moreover,  Armstrong, Cairns and Jessup [ACJ] studied the cohomology
of certain $2$-step nilpotent extensions of abelian Lie algebras.
Furthermore, Armstrong and Sigg [AS] generalized the latter to the
nilpotent Lie algebras which have an abelian ideal with codimension
$1$. It is also a generalization of Bordemann's result [B] on the
cohomology of filiform Lie algebras. All these algebras are
isomorphic to certain tree diagram Lie algebras, and the
 methods of computation used in all of the above papers are
essentially based on the Dixmier's sequences. These sequences are
also used by Fialowski and Millionschikov [FM] in dealing with the
cohomology  of the graded Lie algebras of maximal class, which are
infinite dimensional nilpotent Lie algebras.

Bott [Br] showed that the Betti numbers of the maximal nilpotent
subalgebras of finite dimensional simple Lie algebras can be
expressed by means of the Weyl group. He also pointed out that his
theorem is equivalent to the Weyl  denominator identity. There are
several distinct proofs on Bott's result. For example, it was proved
in [BBG] by representation theory and in [K] by Hodge Laplacian.
Both of these two methods are very effective. The calculation in
[BBG] is generalized by Garland and Lepowsky to the case of
Kac-Moody algebras and the celebrated Macdonald identities were
recovered. The main tool in this paper is the Hodge Laplacian
introduced by Kostant [K].

Euler-Poincar$\acute{\mbox{e}}$ Principle says that for a Lie
algebra $\mathcal {G}=\bigoplus_{\alpha\in\Gamma}\mathcal
{G}_\alpha$ graded by an additive semigroup $\Gamma$,
\begin{eqnarray*}
\prod_{\alpha\in\Gamma}(1-e^\alpha)^{\dim \mathcal
{G}_\alpha}=\sum_{k=0}^\infty\sum_{\alpha\in\Gamma}(-1)^k\dim
H^k_\alpha(\mathcal {G}) e^\alpha,
\end{eqnarray*}
where $e^\alpha$ are the base elements of the semigroup algebra
$\mathbb{C}[\Gamma]$ (e.g., cf. [KK]). We will use it to obtain our
combinatorial identities. Let $\mathcal {G}$ be a finite dimensional
nilpotent Lie algebra. Dixmier [D] proved  that all Betti numbers of
$\mathcal {G}$ are at least two except the zeroth and the highest
which are one. So there is a lower bound of total rank, i.e. $\dim
H(\mathcal {G})\geq 2\dim\mathcal {G}$. Later, Deninger and Singhof
[DS] showed that the length of a polynomial $P(\mathcal {G})$ gives
a lower bound for $\dim H(\mathcal {G})$. Moreover, there is a
``total rank conjecture" (c.f. [CJP]) which has been open for many
years: $$\dim H(\mathcal {G})\geq2^{\dim C(\mathcal
{G})},\quad\mbox{where $C(\mathcal {G})$ is the centre of $\mathcal
{G}$}.$$ Another conjecture that can be found in literatures is the
{``$b_2$-conjecture"} (c.f. [CJP]):
$$b_2>b_1^2/4\quad\mbox{if $\dim \mathcal {G}>2$},
\quad\mbox{where $b_i=\dim H^i(\mathcal {G})$ are the Betti
numbers}.$$ In this present paper, we will prove these two
conjectures for oriented tree diagram Lie algebras.

Oriented tree diagram Lie algebras are introduced by Xu [X] in order
to study certain evolution partial differential equations. They
provide a new realization of some familiar nilpotent Lie algebras
such as the ones mentioned in the second paragraph.

An {\em oriented tree} is a connected oriented graph without cycles.
It can be described as an ordered pair $\mathcal {T}=(\mathcal
{N},\mathcal {E})$, where $$\mathcal
{N}=\{\iota_1,\iota_2,\ldots,\iota_n\}$$ and $$\mathcal
{E}\subset\{(\iota_i,\iota_j)\mid 1\leq i<j\leq n\}$$ are two
disjoint sets. The elements of $\mathcal {N}$ are called {\em nodes}
while the elements of $\mathcal {E}$ are called {\em oriented
edges}.

We call $\iota$ the {\em root node} if $\{\iota'\mid
(\iota',\iota)\in \mathcal {E}\}=\emptyset$, and the {\em tip node}
if $\{\iota'\mid (\iota,\iota')\in \mathcal {E}\}=\emptyset$. Denote
$$\Lambda=\mbox{the set of root nodes}$$ and $$\Gamma=\mbox{the set of tip
nodes}.$$ Define an {\em oriented tree diagram} $$\mathcal
{T}^{d}=(\mathcal {N},\mathcal {E},d)$$ to be an oriented tree
$\mathcal {T}=(\mathcal {N},\mathcal {E})$ with a weight map
$d:\mathcal {E}\rightarrow\mathbb{Z}_{+}$ (the set of positive
integers). We identify an oriented tree diagram $\mathcal
{T}^d=({\cal N},{\cal E},d)$ with a graph by depicting a small
circle for each node in ${\cal N}$ and $d[(\iota_{i},\iota_{j})]$
segments connecting $i$th circle to $j$th circle for the edge
$(\iota_i,\iota_j)\in{\cal E}$, where the orientation is always from
left to right. For instance, the following figure

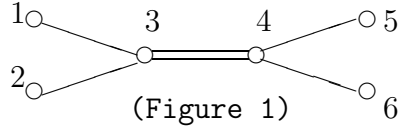
\setlength{\unitlength}{3pt}
\begin{picture}(108,12)
\put(38,4.5){\circle{2}}\put(35,4){$1$}\put(39,4){\line(3,-1){12}}
\put(38,-4.5){\circle{2}}\put(35,-4){$2$}\put(39,-5)
{\line(3,1){12}}
\put(52,0){\circle{2}}\put(52,3){$3$}\put(53,0.5){\line(1,0){12}}
\put(53,-0.5){\line(1,0){12}}
\put(66,0){\circle{2}}\put(66,3){$4$}\put(66.6,0.5)
{\line(3,1){12}}\put(80,4.5){\circle{2}}\put(82,3){$5$}
\put(66.6,-0.5)
{\line(3,-1){12}}\put(80,-4.5){\circle{2}}\put(82,-8){$6$}
\put(50,-8){{\tt(Figure 1)}}
\end{picture}
\vspace{1cm}

\noindent represents the oriented tree diagram $\mathcal
{T}^d=({\cal N},{\cal E},d)$ with ${\cal
N}=\{\iota_1,\iota_2,\ldots,\iota_6\}$, ${\cal
E}=\{(\iota_1,\iota_3),(\iota_2,\iota_3),
(\iota_3,\iota_4),(\iota_4,\iota_5),(\iota_4,\iota_6)\}$, and
$d((\iota_3,\iota_4))=2$, $d({\cal
E}\backslash\{(\iota_3,\iota_4)\})=\{1\}.$

Given a positive integer $n$, there is an associative algebra of
differential operators in $n$ variables:
$$\mathbb{A}=
\sum_{m_{1},m_{2},...,m_{n}=0}^{\infty}
\mathbb{C}[x_{1},x_{2},...,x_
{n}]\partial_{x_{1}}^{m_{1}}\partial_{x_{2}}^{m_{2}}...\partial_{x_{n}}^{m_{n}}.$$
We can define a Lie bracket on $\mathbb{A}$ by
$$[A,B]=AB-BA\quad\mbox{for}\quad \forall A,B\in\mathbb{A}.$$
For any oriented tree diagram
$\mathcal{T}^{d}=(\mathcal{N},\mathcal{E},d)$, we define the Lie
algebra by
$$L_{0}(\mathcal {T}^{d})=\mbox{the
Lie subalgebra of $\mathbb{A}$ generated by }
\{\partial_{x_{i}},x_{j}^{d[(\iota_{j},\iota_{k})]}\partial_{x_{k}}\mid
\iota_i\in \Lambda,\ (\iota_{j},\iota_{k})\in\mathcal {E}\}.$$

Take $\mathcal {T}^{d}$ to be the following diagrams:

\vspace{0.5cm}
\begin{picture}(0,0)
\put(8,6){\circle{2}}\put(8,2){\circle{2}}\put(8,-6){\circle{2}}
\put(4.5,5){$1$}\put(4.5,1){$2$}\put(4.5,-7){$n$}
\put(9,6){\line(2,-1){12}}\put(9,2){\line(6,-1){12}}\put(9,-6){\line(2,1){12}}
\put(10,-3){\vdots}\put(22,0){\circle{2}}\put(20,-3.5){$n+1$}
\put(9,-12){{\tt(Figure 2)}}
\put(46,-4){$n+1$} \put(52,0){\circle{2}}\put(62,-3){\vdots}
\put(53,0){\line(2,1){12}}\put(53,0){\line(6,1){12}}\put(53,0){\line(2,-1){12}}
\put(66,6){\circle{2}}\put(66,2){\circle{2}}\put(66,-6){\circle{2}}
\put(68,5){$1$}\put(68,1){$2$}\put(68,-7){$n$}
\put(50,-12){{\tt(Figure 3)}}
\put(101,0){\circle{2}}\put(106,-0.5){$d$}\put(104,0){.}\put(104,1){.}\put(104,-1){.}
\put(102,-2){\line(1,0){12}}\put(102,2){\line(1,0){12}}\put(102,-1.3){\line(1,0){12}}
\put(100,-5){$1$} \put(115,0){\circle{2}}\put(115,-5){$2$}
\put(98,-12){{\tt(Figure 4)}}
\end{picture}

\vspace{1.2cm}

\noindent It is easy to check that their associated algebras
$L_{0}(\mathcal {T}^{d})$ are the algebras in [S](Heisenberg Lie
algebras), [ACJ] and [B], respectively. When we take $\mathcal
{T}^{d}$ to be the following diagram:

\begin{picture}(108,6)
\put(21,0){\circle{2}}\put(21,-5){1}\put(22,0){\line(1,0){12}}
\put(36,0){\circle{2}}\put(36,-5){2}
\put(42,-0.5){$\cdots$}\put(52,0){\circle{2}}\put(48,-5){$i-1$}\put(53,0){\line(1,0){12}}
\put(66,0){\circle{2}}\put(66,-5){$i$} \put(35,-12){{\tt(Figure 5)}}
\end{picture}
\vspace{1.5cm}

\noindent the Lie algebra $L_0(\mathcal {T}^d)$ is just the maximal
nilpotent subalgebras of $sl(i+1)$.

A natural generalization of both Figure 3 and Figure 5 is the
diagram $\mathcal {A}_n^m$:

\vspace{0.5cm}
\begin{picture}(108,6)
\put(21,0){\circle{2}}\put(21,-5){1}\put(22,0){\line(1,0){12}}
\put(36,0){\circle{2}}\put(36,-5){2}
\put(42,-0.5){$\cdots$}\put(52,0){\circle{2}}\put(48,-5){$n-1$}\put(53,0){\line(1,0){12}}
\put(66,0){\circle{2}}\put(66,-5){$n$} \put(44,-12){{\tt(Figure 6)}}
\put(66,0){\circle{2}}\put(76,-3){\vdots}
\put(67,0){\line(2,1){12}}\put(67,0){\line(6,1){12}}\put(67,0){\line(2,-1){12}}
\put(80,6){\circle{2}}\put(80,2){\circle{2}}\put(80,-6){\circle{2}}
\put(82,5){$n+1$}\put(82,1){$n+2$}\put(82,-7){$n+m$}
\end{picture}
\vspace{1.5cm}

\noindent Without confusion, we also identify $L_0(\mathcal
{A}_n^m)$ with $\mathcal {A}_n^m$ for short. In this paper, We will
compute $H(\mathcal {A}_n^m)$. The Betti numbers of $\mathcal
{A}_1^m$ had been obtained in [ACJ], which is a very special case of
ours.

The paper is organized as follows. In Sections 2 and 3, we review
all necessary definitions and the known facts concerning oriented
tree diagram Lie algebras $L_0(\mathcal {T}^d)$ and their
cohomology, especially the Hodge Laplacian introduced  by Kostant
[K]. We also use the Hodge Laplacian to prove that both of the total
rank conjecture and $b_2$-conjecture hold for any oriented tree
diagram Lie algebras $L_0(\mathcal {T}^d)$ at the end of Section 3.
In Section 4, $H(\mathcal {A}_n^m)$ is computed and then an analogue
of the Weyl denominator identity is obtained by
Euler-Poincar$\acute{\mbox{e}}$ principle, where the Vandermonde
determinant identity is a special case. The last section is devoted
to the calculation of the cohomology the solvable Lie algebra
$L_1(\mathcal
{T}^d)=\sum_{i=1}^n\mathbb{C}x_i\partial_{x_i}+L_0(\mathcal {T}^d)$.

\section{Notations and Facts on $L_0(\mathcal {T}^d)$}

Given an oriented tree diagram $\mathcal {T}^{d}=(\mathcal
{N},\mathcal {E},d)$, for $\forall\iota_i,\iota_j\in \mathcal {N}$,
denote $$\mathcal
{C}_{i,j}=\{\iota_{i_1}=\iota_i,\iota_{i_2},\ldots,\iota_{i_r}=\iota_j\}$$
to be the sequence of nodes with
$$(\iota_{i_{1}},\iota_{i_{2}}),(\iota_{i_{2}},\iota_{i_{3}}),...,
(\iota_{i_{r-1}},\iota_{i_{r}})\in\mathcal {E}.$$ We remake that
$\mathcal {C}_{i,j}$ is unique determined by $\iota_i$ and
$\iota_j$. Of course, sometimes $\mathcal {C}_{i,j}$ may be
$\emptyset$. We denote $\mathcal {C}_{i,i}=\{\iota_i\}$ for
convenience.

Set $$\mathcal {C}_i=\{\iota_j\mid \mathcal {C}_{j,i}\neq
\emptyset\},\quad\quad \mathcal {D}_i=\{\iota_j\mid
\mathcal{C}_{i,j}\neq \emptyset\},$$ and denote $$\mathcal
{E}_i=\{(\iota_r,\iota_s)\in\mathcal {E}\mid
\iota_r,\iota_s\in\mathcal{C}_i\},\quad\quad \mathcal
{E}_{i,j}=\{(\iota_r,\iota_s)\in\mathcal {E}\mid
\iota_r,\iota_s\in\mathcal{C}_{i,j}\}.$$ Let
$$\kappa_i=\prod_{(\iota_r,\iota_s)\in\mathcal {E}
i}d[(\iota_r,\iota_s)],\quad\quad
\kappa_{i,j}=\frac{\kappa_i}{\prod_{(\iota_r,\iota_s)\in\mathcal {E}
{i,j}}d[(\iota_r,\iota_s)]}.$$ It is obvious that $\mathcal {E}
{i,i}=\emptyset$ and
$\kappa_{i,i}=\kappa_j=1(\forall\iota_i\in\mathcal
{N},\forall\iota_j\in\Lambda)$. Recall that $\Lambda$ is the set of
root nodes and $\Gamma$ is the set of tip nodes. We have a basis of
$L_0(\mathcal {T}^d)$:
$$B(\mathcal {T}^d)=\{\partial_{x_i},(\prod_{\iota_s\in \mathcal
{C}_j\backslash\{\iota_j\}}x_s^{m_s})\partial_{x_j}\mid
\iota_i\in\Lambda, m_s\in\mathbb{N},\sum_{\iota_s\in
\mathcal{C}_j\backslash\{\iota_j\}}m_s\kappa_{s,j}\leq
\kappa_{j}\},$$ where $\mathbb{N}$ is the set of nonnegative
integers. We call it the {\em natural basis} of $L_0(\mathcal
{T}^d)$.

\begin{rem}
For a finite dimensional semisimple Lie algebra $\mathcal {G}$, its
Chevalley basis $B$ possesses a good property: for any $u_1,u_2\in
B$, we always have $[u_1,u_2]=\alpha u_3$, where $u_3\in B$ and
$\alpha\in \mathbb{Z}$. Now we can check easily that the natural
basis $B(\mathcal {T}^d)$ also have this property. In the next
section, this property will help us introduce the Hodge Laplacian.
\end{rem}

The following lemma is obvious and will be used later.
\begin{lem}
{\em\bf ([X], [L])} $\mbox{The center of }L_0({\mathcal
{T}}^d)=\sum_{\iota_i\in \Gamma}\mathbb{C}\partial_{x_i}.$
\hfill$\Box$
\end{lem}

\vspace{0.3cm} In order to describe the result in latter sections
laconically, we add some notations and definitions here.

An oriented tree diagram $\mathcal {T'}^{d'}=(\mathcal {N}',\mathcal
{E}',d')$ is called a {\em subdiagram} of the oriented tree diagram
$\mathcal {T}^d=(\mathcal {N},\mathcal {E},d)$ if $\mathcal
{N}'\subset \mathcal {N}, \mathcal {E}'\subset\mathcal {E}, d'=d|
{\mathcal {E}'}$. Further, if $\mathcal{C}_i\subset\mathcal {N}'$
for any $\iota_i\in \mathcal {N}'\subset\mathcal {N}$ , we call
$\mathcal {T'}^{d'}$ a {\em homo-clan subdiagram} of ${\mathcal
{T}}^{d}$ and call $L_0({\mathcal {T}'}^{d'})$ a {\em homo-clan
subalgebra} of $L_0({\mathcal {T}}^{d})$. For example, Figure 5 is a
homo-clan subdiagram of Figure 6 if $i\leq n+1$.

The following lemma can be got immediately by the definition of
homo-clan subdiagram.
\begin{lem}
If $L_0({\mathcal {T}'}^{d'})$ is a homo-clan subalgebra of
$L_0({\mathcal {T}}^{d})$, then $B(L_0({\mathcal {T}'}^{d'}))\in
B(L_0({\mathcal {T}}^{d}))$. Furthermore, for any $u_1,u_2\in
B(L_0({\mathcal {T}}^{d}))$ with $0\neq[u_1,u_2]\in L_0({\mathcal
{T}'}^{d'})$, we have $u_1,u_2\in B(L_0({\mathcal {T}'}^{d'}))$.
\hfill$\Box$
\end{lem}

Furthermore, there is a graded structure in $L_0({\mathcal
{T}}^{d})$ with $\partial_{x_j}\in L_0({\mathcal
{T}}^{d})_{\epsilon_j-\epsilon_0}$ and
$$(\prod_{\iota_s\in \mathcal
{C} j\backslash\{\iota_j\}}x_s^{m_s})\partial_{x_j}\in L_0({\mathcal
{T}}^{d})_{\epsilon_j-\sum_{\iota_s\in
\mathcal{C}_j\backslash\{\iota_j\}}m_s\epsilon_s}.$$

\vspace{0.3cm} \noindent {\bf Example 1:} Denote $x_0=1$ and
$y_i=x_{n+i}$ for convenience. The Lie algebra $\mathcal {A}_n^m$
(associated with Figure 6) is generated by
$$\{x_{i-1}\partial_{x_i},x_n\partial_{y_j}\mid 1\leq
i\leq n,1\leq j\leq m\}.$$ The natural basis of $\mathcal {A}_n^m$
is
$$B(\mathcal
{A}_n^m)=\{x_{i_1}\partial_{x_{i_2}},x_j\partial_{y_k},\mid 0\leq
i_1<i_2\leq n, 0\leq j\leq n, 1\leq k\leq m\}.$$ Obviously,
$$\mathcal {A}_n^m=\bigoplus_{0\leq i<n;\ i<j\leq m+n}(\mathcal
{A}_n^m)_{\epsilon_j-\epsilon_i},$$ where $(\mathcal
{A}_n^m)_{\epsilon_j-\epsilon_i}=\mathbb{C}x_i\partial x_j$. Hence
$\dim(\mathcal {A}_n^m)_{\epsilon_j-\epsilon_i}=1$.

Set $$A_1=\{\partial_{x_1}\},
A_2=\{\partial_{x_2},x_1\partial_{x_2}\}, \ldots,
A_n=\{\partial_{x_n},x_1\partial_{x_n},\ldots,x_{n-1}\partial_{x_n}\},$$
and $B_{m,n}=\{x_i\partial_{y_ j}\mid 0\leq i\leq n, 1\leq j\leq
m\}.$ We have $B(\mathcal {A}_n^m)=(\bigcup_{i=1}^nA_i)\bigcup
B_{m,n}$.

Let $\mathcal {A}_i^0$ be the algebra associated with Figure 5. By
definition, $\bigcup_{j=1}^iA_j$ is exactly the natural basis of
$\mathcal {A}_i^0$, i.e. $\bigcup_{j=1}^iA_j=B(\mathcal {A}_i^0)$.
Moreover, $\mathcal {A}_i^0\ (0<i\leq n)$ is a homo-clan subalgebra
of both $\mathcal {A}_n^m$ and $\mathcal {A}_k^0\ (k\geq i)$. It is
obvious that $\mathcal {A}_n^0\cong\mathcal {A}_{n-1}^1$.
Furthermore, one can check easily that $\mathcal {A}_n^0$ is
isomorphic to the maximal nilpotent subalgebra of $sl(n+1)$.

\section{Lie Algebra Cohomology and Hodge Laplacian}

Let $\mathcal {G}$ be a finite dimensional Lie algebra over
$\mathbb{F}$ and let $\mathcal {G}^*$ be the vector space dual of
$\mathcal {G}$. The spaces
\\$\wedge \mathcal {G}=\oplus_{i\geq0}\wedge^i\mathcal {G}$ and
$\wedge \mathcal {G}^*=\oplus_{i\geq0}\wedge^i\mathcal {G}^*$ are
their exterior algebras. We have a cochain complex:
$$\mathbb{F}\quad^{\underrightarrow{D_0}}\ \mathcal {G}^*
\quad ^{\underrightarrow{D_1}}\wedge^2\mathcal {G}^*\quad
^{\underrightarrow{D_2}}\quad\cdots\quad
^{\underrightarrow{D_{i-1}}}\wedge^i\mathcal {G}^*\quad
^{\underrightarrow{D_i}}\cdots.$$ The coboundary operator $D_p$ is
defined by
$$D_p f(r_0,r_1,\ldots,r_p)=\sum_{0\leq i<j\leq
p}(-1)^{i+j}f([r_i,r_j],r_0,\ldots,\widehat{r_i},\ldots,\widehat{r_j},\ldots,r_p),$$
where the sign\quad $\widehat{}$\quad indicates that the argument
below it must be omitted.

The cohomology of $(\wedge\mathcal {G}^*,D)$ is called the {\em
cohomology} (with trivial coefficients) of the Lie algebra $\mathcal
{G}$ and is denoted by $H(\mathcal {G})$. The gradation from
$\mathcal {G}$ induces a gradation in $H(\mathcal {G})$. Given a
basis $B=\{u_1,u_2,...,u_n\}$ of $\mathcal {G}$, denote
$\{u_1^*,u_2^*,...,u_n^*\}$ the basis of $\mathcal {G}^*$ where
$u_i^*$ is the linear function on $\mathcal {G}$ with
$u_i^*(u_j)=\delta_{i,j}$. In rest of this paper, we always denote
by $\mathbb{F}\langle X\rangle$ the  polynomial algebra generated by
fermionic variables in $X$ with operation ``$\wedge$''. Now we
suppose $\mathcal {G}=L_0(\mathcal {T}^d)$ and take
$B=\{u_1,u_2,...,u_n\}$ to be its natural basis $B(\mathcal {T}^d)$.
Then $\wedge L_0(\mathcal {T}^d)=\mathbb{F}\langle B(L_0(\mathcal
{T}^d))\rangle$.

Recall the property mentioned in Remark 2.1: for $\forall u_i,u_j\in
B$, $[u_i,u_j]=\alpha u_k$, where $u_k\in B$ and
$\alpha\in\mathbb{Z}$. By this property and under the nature
isomorphism $\wedge^i\mathcal {G}^*\cong(\wedge^i\mathcal {G})^*$,
it is easy to check that $D_p$ is the linear map with
$$D_p(u_{i_1}^*\wedge u_{i_2}^*\wedge\cdots\wedge u_{i_t}^*)
=\sum_{k=1}^t(-1)^k u_{i_1}^*\wedge
u_{i_2}^*\wedge\cdots\wedge(\Delta u_{i_k})^*\wedge\cdots\wedge
u_{i_t}^*,$$ where\ $\Delta:\mathcal {G}\rightarrow\wedge^2\mathcal
{G}$ is the linear map with
\begin{equation}
\Delta(u_i)=\sum_{[u_j,u_k]=\alpha u_i,\;\alpha\in\mathbb{Z}
+}\alpha u_j\wedge u_k.
\end{equation}
Without confusion, we can identify $\wedge^i\mathcal {G}^*$ with
$\wedge^i\mathcal {G}$ and redefine
$$D_p(r_1\wedge r_2\wedge\cdots\wedge r_p)
=\sum_{i=1}^p(-1)^i r_1\wedge r_2\wedge\cdots\wedge(\Delta
r_i)\wedge\cdots\wedge r_p,$$ i.e. $D_p:\wedge^p\mathcal
{G}\rightarrow\wedge^{p+1}\mathcal {G}$ is the coboundary operator
of complex:
$$\mathbb{F}\quad^{\underrightarrow{D_0}}\ \mathcal {G}
\quad ^{\underrightarrow{D_1}}\wedge^2\mathcal {G}\quad
^{\underrightarrow{D_2}}\quad\cdots\quad
^{\underrightarrow{D_{i-1}}}\wedge^i\mathcal {G}\quad
^{\underrightarrow{D_i}}\cdots.$$ On the other hand, there also
exists a chain complex:
$$\mathbb{F}\quad^{\underleftarrow{\delta_1}}\ \mathcal {G} \quad
^{\underleftarrow{\delta_2}}\wedge^2\mathcal {G}\quad
^{\underleftarrow{\delta_3}}\quad\cdots\quad
^{\underleftarrow{\delta_{i}}}\wedge^i\mathcal {G}\quad
^{\underleftarrow{\delta_{i+1}}}\cdots,$$ where the boundary
operator $\delta$ is defined by
$$\delta_p(r_1\wedge r_2\wedge\cdots\wedge r_p)=
\sum_{1\leq i<j\leq p}(-1)^{i+j} [r_i,r_j]\wedge
r_1\wedge\cdots\wedge\widehat{r_i}\wedge\cdots\wedge\widehat{r_j}\wedge\cdots\wedge
r_p.$$

Now we can define the operator $$\mathcal {L}=D\delta+\delta
D:\wedge\mathcal {G}\rightarrow\wedge\mathcal {G}$$ and call it the
{\em Hodge Laplacian}. Precisely, $$\mathcal {L}_p=\mathcal
{L}\mid_{\wedge^p\mathcal
{G}}=D_p\delta_{p+1}+\delta_pD_{p-1}:\wedge^p\mathcal
{G}\rightarrow\wedge^p\mathcal {G}.$$

\begin{thm}
{\em\bf([K])} One has a direct sum (a ``Hodge decomposition''),
$$\wedge \mathcal {G}=\mbox{Im}\ \mathcal {L}\oplus\mbox{Ker}\ \mathcal
{L},\quad\mbox{(hence $\wedge^p\mathcal {G}=\mbox{Im}\ \mathcal
{L}_{p}\oplus\mbox{Ker}\ \mathcal {L}_p$)}$$ and $$\mbox{Ker}\
\mathcal {L}_p=\mbox{Ker}\ D_p\cap\mbox{Ker}\ \delta_p\ ;\quad
\mbox{Im}\ \mathcal {L}_p=\mbox{Im}\ D_{p-1}\oplus\mbox{Im}\
\delta_{p+1}.$$ \hfill$\Box$
\end{thm}

Elements in $\mbox{Ker}\ \mathcal {L}$ are called {\em harmonic}. By
the above theorem, $c\in \wedge\mathcal {G}$ is harmonic if and only
if $D(c)=0,\ \delta(c)=0$. For convenience, we also denote
$\widetilde{H}^p(\mathcal {G})=\mbox{Ker}\ \mathcal {L}_p$ and
$\widetilde{H}(\mathcal {G})=\mbox{Ker}\ \mathcal {L}$.

\begin{thm}
{\em\bf([K],[F])} Every element of the space $H^p(\mathcal {G})$ can
be represented by a unique harmonic cocycle from $\wedge\mathcal
{G}$, namely, there is a natural isomorphism
$$\widetilde{H}^p(\mathcal {G})=\mbox{Ker}\ \mathcal {L}_p\cong
H^p(\mathcal {G}).$$ \hfill$\Box$
\end{thm}

By the above two theorems, we have:
\begin{lem}
If $\mathcal {T}'^{d'}$ is a homo-clan subdiagram of $\mathcal
{T}^d$, then $\widetilde{H}(L_0(\mathcal {T}'^{d'}))\subset
\widetilde{H}(L_0(\mathcal {T}^{d}))$.
\end{lem}
\begin{proof}
By Lemma 2.3, it follows from the fact that the coboundary operator
$D$ and the boundary operator $\delta$ of $L_0(\mathcal {T}'^{d'})$
are just the ones of $L_0(\mathcal {T}^{d})$ restricted to
$L_0(\mathcal {T}'^{d'})$.
\end{proof}

\begin{thm} The total rank conjecture
$\dim H(L_{0}(\mathcal {T}^{d}))\geq2^{\dim C(L_{0}(\mathcal
{T}^{d}))}$ holds.
\end{thm}
\begin{proof}
By Lemma 2.1, $\dim C(L_{0}(\mathcal {T}^{d}))=|\Gamma|$ (i.e. the
number of elements in $\Gamma$). For the trivial case $|\mathcal
{N}|=1$, it is true that $\dim H(L_{0}(\mathcal
{T}^{d}))=2\geq2^{\dim C(L_{0}(\mathcal {T}^{d}))}=2$. Thus we
assume $|\mathcal {N}|>1$, and hence
$\Lambda\bigcap\Gamma=\emptyset$.

Suppose $\Gamma=\{\iota_1,\iota_2,\ldots,\iota_t\}$, then
$|\Gamma|=t$. For any $\iota_i\in\Gamma(i=1,2,\ldots,t)$, we can
take a $\iota_{p(i)}\in\mathcal {N}$ such that
$d[(\iota_{p(i)},\iota_i)]=d_i\neq0$. Take
$A=\{x_{p(i)}^{D_i}\partial_{x i}\mid i=1,2,\ldots,t\}\subset
B(L_{0}(\mathcal {T}^{d}))$ and $\wedge A$ its exterior algebra.

Recall the operator $\Delta:\mathcal {G}\rightarrow\wedge^2\mathcal
{G}$ introduced in (3.1). As $A\bigcap[L_{0}(\mathcal
{T}^{d}),L_{0}(\mathcal {T}^{d})]=\emptyset$, we have
$\Delta(x_{p(i)}^{d_i}\partial_{x_i})=0$. Thus $D(\wedge A)=0$. On
the other hand, $[A,A]=\{0\}$, so $\delta(\wedge A)=0$. Thanks to
Theorem 3.1 and 3.2, we have $\wedge A\subset\mbox{Ker} \mathcal
{L}\cong H(L_{0}(\mathcal {T}^{d}))$. So $\dim H(L_{0}(\mathcal
{T}^{d}))\geq\dim (\wedge A)=2^t=2^{\dim C(L_{0}(\mathcal
{T}^{d}))}$.
\end{proof}

\vspace{0.6cm} Next we turn to the $b_2$-conjecture.
\begin{thm} For any $\mathcal {T}^{d}$ except the trivial case $|\mathcal
{N}|=1$, we always have $b_2>b_1^2/4$, where $b_i=\dim
H^i(L_{0}(\mathcal {T}^{d}))$ are the Betti numbers.
\end{thm}
\begin{proof}
It is obvious that $\mbox{Ker}\:D_1\bigcap\mbox{Ker}\: \delta_1
=\mbox{Span} \{\partial_{x_i},x_j^{d_[(\iota_j,\iota_k)]}\partial
{x_k}\mid\iota_i\in\Lambda,(\iota_j,\iota_k)\in\mathcal {E}\}$. Thus
by Theorem 3.1 and 3.2, $b_1=\dim H^1(L_{0}(\mathcal
{T}^{d}))=|\Lambda|+|\mathcal {E}|$.

On the other hand, one can check that
\begin{eqnarray*}
\{\partial_{x_{i_1}}\wedge\partial_{x_{i_2}},\partial_{x_{i}}\wedge
x_{j}^{d_[(\iota_{j},\iota_{k})]}\partial_{x_{k}},
x_{j_1}^{d_[(\iota_{j_1},\iota_{k_1})]}\partial_{x_{k_1}}\wedge
x_{j_ 2}^{d_[(\iota_{j_2},\iota_{k_2})]}\partial_{x_{k_2}}\mid
\iota_i,\iota_{i_1},\iota_{i_2}\in\Lambda,\quad\quad\quad\quad\quad\quad\\
(\iota_{j},\iota_{k}),(\iota_{j_1},\iota_{k_1}),(\iota_{j_1},\iota_{k_1})\in\mathcal
{E},i\neq j,j_1\neq k_2,j_2\neq k_1\}\\\subset
\mbox{Ker}\:D_2\bigcap\mbox{Ker}\:\delta_2=\mbox{Ker}\:\mathcal
{L}_2,\quad\quad\quad\quad\quad\quad\quad\quad\quad\quad\quad
\quad\quad\quad\quad\quad\quad\quad\quad\quad\quad\quad\quad\quad\quad\quad
\end{eqnarray*}
and
\begin{eqnarray*}
\{x_{j}^{d_[(\iota_{j},\iota_{k})]-1}\partial_{x_{k}}\wedge
x_{j}^{d_[(\iota_{j},\iota_{k})]}\partial_{x_{k}}\mid
(\iota_{j},\iota_{k})\in\mathcal {E}\}\subset
\mbox{Ker}\:D_2\bigcap\mbox{Ker}\:\delta_2=\mbox{Ker}\:\mathcal
{L}_2.
\end{eqnarray*}
Hence $b_2=\dim H^2(L_{0}(\mathcal {T}^{d}))\geq{|\Lambda|+|\mathcal
{E}|\choose 2}-|\mathcal {E}|+|\mathcal {E}|={|\Lambda|+|\mathcal
{E}|\choose 2}.$

So if $|\Lambda|+|\mathcal {E}|>2$, it must be $b_2>b_1^2/4$.
Otherwise, when $|\Lambda|+|\mathcal {E}|=2$, $\mathcal {T}^d$ must
be Figure 4. It is easy to get $$\mbox{Ker}\:\mathcal
{L}_1=\mbox{Span}\{\partial_{x_1},x_1^{d_[(\iota_1,\iota_2)]}\partial_{x_2}\}$$
and
$$\{\partial_{x_1}\wedge\partial_{x_2},
x_1^{d_[(\iota_1,\iota_2)]-1}\partial_{x_2}\wedge
x_1^{d_[(\iota_1,\iota_2)]}\partial_{x_2}\}
\subset\mbox{Ker}\:\mathcal {L}_2,$$ i.e. $b_1=2$ and $b_2\geq2$.
There also have to be $b_2>b_1^2/4$.
\end{proof}

\section{Cohomology of $\mathcal {A}_n^m$}

In this section, we will compute the cohomology of $\mathcal
{A}_n^m$. In other words, we want to get all the harmonic cocycle of
$\wedge\mathcal {A}_n^m$ due to the Theorem 3.2.

With the notations introduced to $\mathcal{A}_n^m$ in Example 1 at
the end of Section 2, we define an ordering ``$\prec$" on
$B(\mathcal{A}_n^m)$ by
$$x_i\partial_{x_j}\prec x_k\partial_{x_l}\quad\mbox{if}\quad j<l\ \mbox{or}\ j=l\ \mbox{and}\ i<k.$$

Set $$\mathscr{B}(\mathcal{A}_n^m)=\{u_{1}\wedge
u_{2}\wedge\cdots\wedge u_{p}\mid u_{1}\prec u_{2}\prec\cdots\prec
u_p, \;p=0,1,2,\ldots\}.$$ Obviously, $\mathscr{B}(\mathcal{A}_n^m)$
forms  a basis of $\wedge \mathcal{A}_n^m$, which can be regarded as
the polynomial algebra $\mathbb{F}\langle B(\mathcal{A}_n^m)\rangle$
of fermionic variables.

\begin{lem}
For any $a\wedge b\in \wedge \mathcal {A}_n^m$ with $0\neq
a\in\mathbb{F}\langle B(\mathcal{A}_n^0)\rangle$ and $0\neq
b=x_{i_1}\partial_{y_{j_1}}\wedge\cdots\wedge
x_{i_k}\partial_{y_{j_k}}$, we have $\mathcal {L}(a\wedge b)=a\wedge
b'$ with $b'\in\mathbb{F}\langle B_{m,n}\rangle$. Moreover, each
term of $b'$ must be of the form
$\alpha(x_{i_{\sigma_{s,t}(1)}}\partial_{y_{j_1}}\wedge\cdots\wedge
x_{i_{\sigma_{s,t}(k)}}\partial_{y_{j_k}})$, where $\alpha\in
\mathbb{F}$ and $\sigma_{s,t}(s\leq t) \in S_k$ satisfy
$\sigma_{s,t}(s)=t$, $\sigma_{s,t}(t)=s$ and $\sigma_{s,t}(r)=r
(r\neq s,r\neq t)$. In particular, all of the monomials in
$\mathbb{F}\langle B(\mathcal{A}_n^0)\rangle$ are eigenvectors of
the linear map $\mathcal {L}$. Thus $\mbox{Ker}\mathcal\:
{L}\bigcap\mathscr{B}(\mathcal{A}_n^0)$ forms a basis of
$\widetilde{H}(\mathcal{A}_n^0)$.
\end{lem}
\begin{proof}
By the definitions of $D$ and $\delta$, there are at most three
distinct factors between  any term of $D\delta (a\wedge b)$ and any
term of $a\wedge b$. So are there between any term of $\delta
D(a\wedge b)$ and  any term of $a\wedge b$.

For any term of $D\delta (a\wedge b)$  with three distinct factors
relative to $a\wedge b$ , it must be produced by
$$x_i\partial_{x_j}\wedge x_j\partial_{z_1}\wedge
x_s\partial_{z_2}\wedge\cdots \stackrel{\delta}{\dashrightarrow}
-x_i\partial_{z_1}\wedge
x_s\partial_{z_2}\wedge\cdots\stackrel{D}{\dashrightarrow}
-x_i\partial_{z_1}\wedge x_s\partial_{x_u}\wedge
x_u\partial_{z_2}\wedge\cdots$$ where $z_1,z_2$ may be $x_k$ or
$y_k$.

(In this paper, the ``$\dashrightarrow$'' but not the
``$\rightarrow$'' will be used in the calculation frequently. If we
write ``$a \stackrel{f}{\dashrightarrow} b$'', $b$ may be not equal
to $f(a)$ but equal to the terms of $f(a)$ which we are concerned
about.)

At the same time, there must be a term of $\delta D(a\wedge b)$
produced by
$$x_i\partial_{x_j}\wedge x_j\partial_{z_1}\wedge
x_s\partial_{z_2}\wedge\cdots\stackrel{D}{\dashrightarrow}
-x_i\partial_{x_j}\wedge x_j\partial_{z_1}\wedge
x_s\partial_{x_u}\wedge
x_u\partial_{z_2}\wedge\cdots\stackrel{\delta}{\dashrightarrow}x_i\partial_{z_1}\wedge
x_s\partial_{x_u}\wedge x_u\partial_{z_2}\wedge\cdots.$$ These two
terms counteract each other. For any term of $D\delta(a\wedge b)$
with two distinct factors relative to $a\wedge b$ , it may be
produced by
$$x_i\partial_{x_j}\wedge x_j\partial_{z}\wedge\cdots \stackrel{\delta}{\dashrightarrow}
-x_i\partial_{z}\wedge\cdots\stackrel{D}{\dashrightarrow}
x_i\partial_{x_s}\wedge x_s\partial_{z}\wedge\cdots$$ where $z$ may
be $x_k$ or $y_k$.

Suppose $s<j$ (the case of $s>j$ can be checked similarly). If
$x_s\partial_{x_j}$ is also a factor of $a\wedge b$, then there must
be another term of $D\delta (a\wedge b)$ produced by
$$x_i\partial_{x_j}\wedge x_j\partial_{z}\wedge x_s\partial_{x_j}\wedge\cdots\stackrel{\delta}{\dashrightarrow}
-x_i\partial_{x_j}\wedge
x_s\partial_{z}\wedge\cdots\stackrel{D}{\dashrightarrow}
-x_i\partial_{x_s}\wedge x_s\partial_{z}\wedge
x_s\partial_{x_j}\wedge\cdots.$$ If $x_s\partial_{x_j}$ is not a
factor of $a\wedge b$, then there must be a term of $\delta
D(a\wedge b)$ produced by
$$x_i\partial_{x_j}\wedge x_j\partial_{z}\wedge\cdots \stackrel{D}{\dashrightarrow}
-x_i\partial_{x_s}\wedge x_s\partial_{x_j}\wedge
x_j\partial_{z}\wedge\cdots\stackrel{\delta}{\dashrightarrow}-x_i\partial_{x_s}\wedge
x_s\partial_{z}\wedge\cdots.$$ Thus these terms can always be
counteracted.

But the following two can not be counteracted. One is produced by:
(when $x_i\partial_{x_j}$ is a factor of $a$)
\begin{equation*}x_i\partial_{x_j}\wedge x_j\partial_{y_k}\wedge
x_i\partial_{y_l}\wedge\cdots
\stackrel{\delta}{\dashrightarrow}-x_i\partial_{y_k}\wedge
x_i\partial_{y_l}\wedge\cdots\stackrel{D}{\dashrightarrow}
x_i\partial_{x_j}\wedge x_i\partial_{y_k}\wedge
x_j\partial_{y_l}\wedge\cdots.\end{equation*} The other is produced
by: (when $x_i\partial_{x_j}$ is not a factor of $a$)
\begin{equation*}x_j\partial_{y_k}\wedge x_i\partial_{y_l}\wedge\cdots
\stackrel{D}{\dashrightarrow}-x_i\partial_{x_j}\wedge
x_j\partial_{y_k}\wedge
x_j\partial_{y_l}\wedge\cdots\stackrel{\delta}{\dashrightarrow}
x_i\partial_{y_k}\wedge x_j\partial_{y_l}\wedge\cdots.
\end{equation*}
Moreover, it is obvious that there can not be a common term of
$D\delta(a\wedge b)$ and $\delta D(a\wedge b)$ that has only one
distinct factor relative to $a\wedge b$. So $\mathcal {L}(a\wedge
b)$ must be of the form described in the lemma.

In particular, if we take $b=1$, then  all of the monomials in
$\mathbb{F}\langle B(\mathcal{A}_n^0)\rangle$ are eigenvectors of
$\mathcal {L}$. As $\mathscr{B}(\mathcal{A}_n^0)$ is a basis of
$\mathbb{F}\langle B(\mathcal{A}_n^0)\rangle$, $\mbox{Ker}\:\mathcal
{L}\bigcap\mathscr{B}(\mathcal{A}_n^0)$ forms a basis of
$\widetilde{H}(\mathcal{A}_n^0)$.
\end{proof}

\begin{lem}
For any $0\neq a\wedge b\in \widetilde{H}(\mathcal {A}_n^m)$ with
$a\in\mathscr{B}(\mathcal{A}_{n}^0)$ and $0\neq b\in\wedge B_{m,n}$,
we have $a\in\widetilde{H}(\mathcal {A}_{n}^0)$. In particular, if
$0\neq a\wedge b\in \widetilde{H}(\mathcal {A}_i^0)$ with
$a\in\mathscr{B}(\mathcal{A}_{i-1}^0)$ and $0\neq b\in\wedge A_i$,
then $a\in\widetilde{H}(\mathcal {A}_{i-1}^0)$.
\end{lem}
\begin{proof}
By Theorem 3.1, we only need to prove that both $D(a)$ and
$\delta(a)$ are zero.

$0=D(a\wedge b)=D(a)\wedge b+(-1)^{\deg(a)}a\wedge D(b)$, where
$\deg(a)$ is the degree of $a$  as a monomial of $\mathbb{F}\langle
B(\mathcal{A}_{i-1}^0)\rangle $. As $a$ is not contained in any term
of $D(a)\wedge b$, we must have $D(a)\wedge b=a\wedge D(b)=0$. Hence
$D(a)=0$.

$0=\delta(a\wedge b)=\delta(a)\wedge b+U$. Here the $U$ is a sum of
terms with form $a'\wedge b'$, where
$a'\in\mathscr{B}(\mathcal{A}_{i-1}^0)$ and $b'\in\wedge A_i$.
Furthermore, we can easily observe that each $a'$ loses at most one
factor of $a$ but each term of $\delta(a)$ loses two. Thus it must
be $\delta(a)\wedge b=U=0$. Hence $\delta(a)=0$.
\end{proof}

\begin{cor}
(1). If $\sum\limits_{j=1}^{t} a_j\wedge
b_j\in\widetilde{H}(\mathcal {A}_n^0)$, where
$a_j\in\mathscr{B}(\mathcal{A}_{n-1}^0)$, $b_j\in\wedge A_n$ and
$a_{j_{1}}\neq a_{j_2} (j_1\neq j_2)$. Then for any
$j\in\{1,2,\ldots,t\}$, we have $a_j\wedge
b_j\in\widetilde{H}(\mathcal {A}_n^0)$ and $a_j\in
\widetilde{H}(\mathcal {A}_{n-1}^0)$.

(2). If $0\neq a_1\wedge a_2\wedge\cdots\wedge a_n\in
\widetilde{H}(\mathcal {A}_n^0)$ where $a_i\in\wedge A_i\
(i=1,2,\ldots,n)$, then $a_1\wedge a_2\wedge\cdots\wedge
a_j\in\widetilde{H}(\mathcal {A}_j^0)\quad (j=1,2,\ldots,n)$.
\end{cor}
\begin{proof}
Both of these two statements can be obtained by Lemmas 4.1 and 4.2
directly.
\end{proof}

\vspace{0.5cm} Thanks to the first claim of the above corollary,  we
only need to consider the monomials in $\mathbb{F}\langle
B(\mathcal{A} i^0)\rangle$  in order to obtain a basis of
$\widetilde{H}(\mathcal {A}_{i}^0)$. Precisely,
$\mathscr{B}(\mathcal{A}_{i}^0)\bigcap \widetilde{H}(\mathcal
{A}_{i}^0)$ is a basis of $\widetilde{H}(\mathcal {A}_{i}^0)$.

For each
$a\in\mathscr{B}(\mathcal{A}_i^0)\bigcap\widetilde{H}(\mathcal
{A}_{i}^0)$, we introduce a total ordering ``$\prec_a$'' into
$\{0,1,2,\ldots,i\}$:
\begin{equation}\mbox{for any}\ 0\leq j<k\leq i,\
\left\{\begin{array}{ll} k\prec_a j, &\mbox{if $x_j\partial_{x_k}$ is a factor of $a$};\\
j\prec_a k, &\mbox{otherwise}.
\end{array}\right.
\end{equation}

Although we have not checked that ``$\prec_a$'' is well defined, it
will be indicated in the next theorem.

\begin{thm}
The total orderings defined in (4.2) are well defined. For any
$a\wedge b\in\mathscr{B}(\mathcal{A}_{i+1}^0)\bigcap
\widetilde{H}(\mathcal {A}_{i+1}^0)$ with $0\neq a\in\wedge\mathcal
{A}_{i}^0$ and $0\neq b\in\wedge A_{i+1}$, if
$x_j\partial_{x_{i+1}}$ is a factor of $b$, then all
$x_k\partial_{x_{i+1}}\ (j\prec_a k)$ are factors of $b$.
Conversely, each element of this form must be in
$\mathscr{B}(\mathcal{A}_{i+1}^0)\bigcap \widetilde{H}(\mathcal
{A}_{i+1}^0)$.
\end{thm}

\begin{proof} We will prove it by induction on $i$.
For $i=1$, there are only two elements $1$ and $\partial_{x_1}$ in
$\mathscr{B}(\mathcal{A}_1^0)\bigcap\widetilde{H}(\mathcal
{A}_{1}^0)$. We have
$$0\prec_1 1;\quad 1\prec_{\partial_{x_1}}0.$$
Thus ``$\prec_1$'' and ``$\prec_{\partial_{x_1}}$'' are indeed total
orderings. For $a\wedge b\in\mathscr{B}(\mathcal{A}_{2}^0)\bigcap
\widetilde{H}(\mathcal {A}_{2}^0)$ with $0\neq a\in\wedge\mathcal
{A}_{1}^0$ and $0\neq b\in\wedge A_{2}$, an easy calculation
indicates that $a\wedge b=x_1\partial_{x_2}$, $\partial_{x_2}\wedge
x_1\partial_{x_2}$, $\partial_{x_1}\wedge\partial_{x_2}$, or
$\partial_{x_1}\wedge\partial_{x_2}\wedge x_1\partial_{x_2}$. One
can check that the statement is true for $i=1$.

Suppose that the statement is true for $i-1$. Take any $a\wedge
b\in\mathscr{B}(\mathcal{A}_{i+1}^0)\bigcap \widetilde{H}(\mathcal
{A}_{i+1}^0)$ with $0\neq a\in\wedge\mathcal {A}_{i}^0$ and $0\neq
b\in\wedge A_{i+1}$. By Lemma 4.2, one has
$a\in\mathscr{B}(\mathcal{A}_{i}^0)\bigcap \widetilde{H}(\mathcal
{A}_{i}^0)$. Assume $a=a_1\wedge a_2$ with $0\neq
a_1\in\wedge\mathcal {A}_{i-1}^0$ and $0\neq a_2\in\wedge A_{i}$. By
definition, we know that ``$\prec_{a}$'' restricted to
$\{1,2,\ldots,i-1\}$ should be equal to ``$\prec_{a_1}$''. We can
assume that $j\in\{1,2,\ldots,i-1\}$ is the unique element such that
$x_j\partial_{x_i}$ is a factor of $a_2$ but $x_k\partial_{x_i}\
(\forall\ k \prec_{a_1}j)$ are not. Hence we get that $j$ must be
the next number of $i$ under the total ordering ``$\prec_a$''.
Precisely, $$i\prec_{a}k\ \mbox{if and only if $k=j$ or
$j\prec_{a_1}k$\ (hence $k\prec_a i$ if and only if
$k\prec_{a_1}j$)}.$$ So we have proved that ``$\prec_a$'' is well
defined.

Taking any pair $(k,l)$ such that $x_k\partial_{x_{i+1}}$ is a
factor of $b$ and $k\prec_al$, we need to prove that
$x_l\partial_{x_{i+1}}$ is a factor of $b$. Suppose not.

If $k<l$, then $x_k\partial_{x_l}$ is not a factor of $a$. Thus
there must be a nonzero term of $D(a\wedge b)$:
$$x_k\partial_{x_{i+1}}\wedge\cdots\stackrel{D}{\dashrightarrow}-x_k\partial_{x_l}\wedge x_l\partial_{x_{i+1}}\wedge\cdots.$$
This  contradicts Theorem 3.1. If $l<k$, then $x_l\partial_{x_k}$ is
a factor of $a$. Thus there must be a nonzero term of
$\delta(a\wedge b)$:
$$x_l\partial_{x_k}\wedge x_k\partial_{x_{i+1}}\wedge\cdots\stackrel{\delta}{\dashrightarrow}-x_l\partial_{x_{i+1}}\wedge\cdots.$$
This again leads a contradiction to Theorem 3.1. Therefore
$x_l\partial_{x_{i+1}}$ is a factor of $b$.

At last, if $a\wedge b$ is an element of the form we mentioned. For
any pair $(k,l)$ such that $x_k\partial_{x_{i+1}}$ is a factor of
$b$, we only need to check the following two cases to prove
$D(a\wedge b)=0$ and $\delta(a\wedge b)=0$ (other cases are so
trivial that the total ordering ``$\prec_a$'' is needless). If $l<k$
and $x_l\partial_{x_k}$ is a factor of $a$ (hence $k\prec_a l$),
then $x_l\partial_{x_{i+1}}$ is a factor of $b$. Thus
$$x_l\partial_{x_k}\wedge x_k\partial_{x_{i+1}}\wedge x_l\partial_{x_{i+1}}\wedge\cdots
\stackrel{\delta}{\dashrightarrow}x_l\partial_{x_{i+1}}\wedge
x_l\partial_{x_{i+1}}\wedge\cdots(=0).$$ If $l>k$ and
$x_k\partial_{x_{l}}$ is not a factor of $a$ (hence $k\prec_a l$).
Then $x_l\partial_{x_{i+1}}$ is a factor of $b$. Thus $$
x_k\partial_{x_{i+1}}\wedge
x_l\partial_{x_{i+1}}\wedge\cdots\stackrel{D}{\dashrightarrow}x_k\partial_{x_l}\wedge
x_l\partial_{x_{i+1}}\wedge x_l\partial_{x_{i+1}}\wedge\cdots(=0).$$
Hence $a\wedge b\in\mathscr{B}(\mathcal{A}_{i+1}^0)\bigcap
\widetilde{H}(\mathcal {A}_{i+1}^0)$.
\end{proof}

\begin{cor}
The generating function of the Betti numbers of $\mathcal {A}_n^0$
is
$$\sum_{i=0}^{\infty}b_it^i=(1+t)(1+t+t^2)(1+t+t^2+t^3)\cdots(1+t+t^2+\cdots+t^n)
=\frac{\prod_{i=1}^n(1-t^{i+1})}{(1-t)^n}.$$
\end{cor}
\begin{proof}
For $n=1$, the statement holds trivially.  Suppose the statement is
true for $n-1$. Then for any
$a\in\mathscr{B}(\mathcal{A}_{n-1}^0)\bigcap \widetilde{H}(\mathcal
{A}_{n-1}^0)$, there are $n$ distinct $b\in\wedge A_{n}$ such that
$a\wedge b\in\mathscr{B}(\mathcal{A}_{n}^0)\bigcap
\widetilde{H}(\mathcal {A}_{n}^0)$.

Furthermore, if $\{i_1,i_2,\ldots,i_n\}$ is a permutation of
$\{0,1,\ldots,n-1\}$ such that $i_{n}\prec_a i_{n-1}\prec_a
\cdots\prec_a i_{1}$, then the $n$ distinct $b$ are $1$,
$x_{i_{1}}\partial_{x_n}$, $x_{i_{1}}\partial_{x_n}\wedge
x_{i_{2}}\partial_{x_n}$, $\ldots$ , $x_{i_{1}}\partial_{x_n}\wedge
x_{i_{2}}\partial_{x_n}\wedge\cdots\wedge x_{i_{n}}\partial_{x_n}$,
respectively. Hence the statement holds for $n$.
\end{proof}

\begin{rem}
Bott's theorem has indicated that $b_i=\dim H_i(\mathcal
{A}_n^0)=|S_{n+1}^{(i)}|$ ($|\cdot|$ means the number of elements in
$S_{n+1}^{(i)}$) where $S_{n+1}^{(i)}$ is the set of elements in
$S_{n+1}$ (the $(n+1)$th symmetric group which can be regarded as a
Weyl group) with length $i$. So our generating function is just the
Poincar$\acute{\mbox{e}}$ polynomial which gives the description of
the length about the elements of Weyl group $S_{n+1}$. (The
definition of Poincar$\acute{\mbox{e}}$ polynomial can be found in
[H].) Moreover, we can get the explicit relationship about Bott's
theorem between the cohomology and the Weyl group in the case of
type $A$. That is ,for each
$a\in\mathscr{B}(\mathcal{A}_n^0)\bigcap\widetilde{H}(\mathcal
{A}_{n}^0)$, the total ordering ``$\prec_a$''defined above is
relative to an element $\sigma_a\in S_{n+1}$ with $i\prec_a j
\Leftrightarrow\sigma_a(i)<\sigma_a(j)$.
\end{rem}

\vspace{0.2cm}In order to describe the result about $\mathcal
{A}_n^m$ laconically, we first introduce a notation. Given a total
ordering ``$\prec_a$'' on $\{0,1,\ldots,n\}$, we assume
$i_{n+1}\prec_a i_{n}\prec_a\cdots\prec_a i_1$, where $\{i_1, i_2,
\ldots, i_{n+1}\}$ is a permutation of $\{0,1,\ldots,n\}$.

Set $$Y=\{(j_1,j_2,\ldots,j_k)\mid 1\leq j_1\leq j_2\leq\cdots\leq
j_k\leq m, 0\leq k\leq n\},$$ in which $(j_1,j_2,\ldots,j_k)=0$ if
$k=0$. For any $(j_1,j_2,\ldots,j_k)\in Y$, denote
\begin{eqnarray*}
\varphi^a_{(j_1,j_2,\ldots,j_k)}=\left\{\begin{array}{ll} 1,
&\mbox{if}\quad k=0;\\
\sum_{\sigma\in
S_k}x_{i_1}\partial_{y_{j_{\sigma(1)}}}\wedge\cdots\wedge
x_{i_k}\partial_{y_{j_{\sigma(k)}}},&\mbox{if}\quad k>0.
\end{array}\right.
\end{eqnarray*}

With this notation, we can easily find numbers of elements which
belong to $\widetilde{H}(\mathcal {A}_n^m)$. Denote $$\mathcal
{Q}=\{a\wedge
\varphi^a_{J_1}\wedge\varphi^a_{J_2}\wedge\cdots\wedge\varphi^a_{J_t}\mid
a\in\mathscr{B}(\mathcal{A}_{n}^0)\bigcap\widetilde{H}(\mathcal
{A}_{n}^0), J_1,J_2,\ldots,J_t\in Y, t\in \mathbb{Z}_+\}.$$ We have
\begin{lem}
$\mathcal {Q}\subset\widetilde{H}(\mathcal {A}_n^m)$.
\end{lem}
\begin{proof}
Given any $\varphi^a_{J}=\sum_{\sigma\in
S_k}x_{i_1}\partial_{y_{j_{\sigma(1)}}}\wedge\cdots\wedge
x_{i_k}\partial_{y_{j_{\sigma(k)}}}$, we take the pair $(i_1,i_2)$
(other pairs can be discussed similarly). Suppose $i_1<i_2$, then
$x_{i_1}\partial_{x_{i_2}}$ must be a factor of $a$. Hence
\begin{eqnarray*}x_{i_1}\partial_{x_{i_2}}\wedge
(x_{i_1}\partial_{y_{j_{\sigma(1)}}}\wedge
x_{i_2}\partial_{y_{j_{\sigma(2)}}}+x_{i_1}\partial_{y_{j_{(\sigma\cdot\sigma_{1,2})(1)}}}\wedge
x_{i_2}\partial_{y_{j_{(\sigma\cdot\sigma_{1,2})(2)}}})\wedge\cdots\\\stackrel{\delta}{\dashrightarrow}
(-x_{i_1}\partial_{y_{j_{\sigma(1)}}}\wedge
x_{i_1}\partial_{y_{j_{\sigma(2)}}}-x_{i_1}\partial_{y_{j_{\sigma(2)}}}\wedge
x_{i_1}\partial_{y_{j_{\sigma(1)}}})\wedge\cdots(=0)
\end{eqnarray*}
and
\begin{eqnarray*}
x_{i_1}\partial_{x_{i_2}}\wedge
(x_{i_1}\partial_{y_{j_{\sigma(1)}}}\wedge
x_{i_2}\partial_{y_{j_{\sigma(2)}}}+x_{i_1}\partial_{y_{j_{(\sigma\cdot\sigma_{1,2})(1)}}}\wedge
x_{i_2}\partial_{y_{j_{(\sigma\cdot\sigma_{1,2})(2)}}})\wedge\cdots\stackrel{D}{\dashrightarrow}
\quad\quad\quad\quad\quad\\
x_{i_1}\partial_{x_{i_2}}\wedge (x_{i_1}\partial_{x_{i_2}}\wedge
x_{i_2}\partial_{y_{j_{\sigma(1)}}}\wedge
x_{i_2}\partial_{y_{j_{\sigma(2)}}}+x_{i_1}\partial_{x_{i_2}}\wedge
x_{i_2}\partial_{y_{j_{(\sigma\cdot\sigma_{1,2})(1)}}}\wedge
x_{i_2}\partial_{y_{j_{(\sigma\cdot\sigma_{1,2})(2)}}})\wedge\cdots\\(=0).
\end{eqnarray*}
Suppose $i_1>i_2$, then $x_{i_2}\partial_{x_{i_1}}$ can not be a
factor of $a$. Thus
\begin{eqnarray*}
(x_{i_1}\partial_{y_{j_{\sigma(1)}}}\wedge
x_{i_2}\partial_{y_{j_{\sigma(2)}}}+x_{i_1}\partial_{y_{j_{(\sigma\cdot\sigma_{1,2})(1)}}}\wedge
x_{i_2}\partial_{y_{j_{(\sigma\cdot\sigma_{1,2})(2)}}})\wedge\cdots\stackrel{D}{\dashrightarrow}\quad\quad\quad\\
(x_{i_1}\partial_{y_{j_{\sigma(1)}}}\wedge
x_{i_2}\partial_{x_1}\wedge
x_{i_1}\partial_{y_{j_{\sigma(2)}}}+x_{i_1}\partial_{y_{j_{\sigma(2)}}}\wedge
x_{i_2}\partial_{x_{i_1}}\wedge
x_{i_1}\partial_{y_{j_{\sigma(1)}}})\wedge\cdots(=0)
\end{eqnarray*}
and
\begin{equation*}
(x_{i_1}\partial_{y_{j_{\sigma(1)}}}\wedge
x_{i_2}\partial_{y_{j_{\sigma(2)}}}+x_{i_1}\partial_{y_{j_{(\sigma\cdot\sigma_{1,2})(1)}}}\wedge
x_{i_2}\partial_{y_{j_{(\sigma\cdot\sigma_{1,2})(2)}}})\wedge\cdots\stackrel{\delta}{\dashrightarrow}0.
\end{equation*}

So we always have $D(a\wedge
\varphi^a_{J_1}\wedge\varphi^a_{J_2}\wedge\cdots\wedge\varphi^a_{J_t})=0$
and $\delta(a\wedge
\varphi^a_{J_1}\wedge\varphi^a_{J_2}\wedge\cdots\wedge\varphi^a_{J_t})=0$.
That is $a\wedge
\varphi^a_{J_1}\wedge\varphi^a_{J_2}\wedge\cdots\wedge\varphi^a_{J_t}\in\widetilde{H}(\mathcal
{A}_n^m)$.
\end{proof}

\vspace{0.3cm} By the above lemma, we know $$\mbox{Span}\mathcal
{Q}\subset\widetilde{H}(\mathcal {A}_n^m).$$ The next theorem will
show that the ``$\subset$'' in the above formula can be changed to
``$=$''. In fact, we are even able to take a proper subset $\mathcal
{P}\subset\mathcal {Q}$ such that $\widetilde{H}(\mathcal
{A}_n^m)=\mbox{Span}\mathcal {P}$. Now we first define a such subset
$\mathcal {P}$ and then give the main theorem.

We call $a\wedge
\varphi^a_{J_1}\wedge\varphi^a_{J_2}\wedge\cdots\wedge\varphi^a_{J_t}\in\mathcal
{Q}$ a {\em basic element} if
$J_s=(j_{1,s},j_{2,s},\ldots,j_{p_s,s})\in Y$ $(s=1,2,\ldots,t)$
satisfy that $p_1\geq p_2\geq\cdots\geq p_t$ and
$j_{q,s_1}<j_{q,s_2}$ (with $s_1<s_2$). The set of all basic
elements in $\mathcal {Q}$ is denoted by $\mathcal {P}$.

\begin{thm}
The set $\mathcal {P}$ of all basic elements in $\mathcal {Q}$ is a
basis of $\widetilde{H}(\mathcal {A}_n^m)$.
\end{thm}

We shall divide the proof of Theorem 4.8 into several lemmas.

By Lemma 4.1 and 4.2, we know $\widetilde{H}(\mathcal {A}_n^m)$ can
be spanned by the elements of the form $a\wedge b$ with
$a\in\mathscr{B}(\mathcal{A}_{n}^0)\bigcap\widetilde{H}(\mathcal
{A}_{n}^0)$ and $b\in\wedge B_{m,n}$. Lemma 4.1 also allows us to
fetch $b$ better. Precisely, we can take $b$ to satisfy the
following condition:

\vspace{0.2cm} \noindent{\bf Condition ($\ast$)} {\em If
$x_{s_1}\partial_{y_{t_1}}\wedge\cdots\wedge
x_{s_k}\partial_{y_{t_k}}$ is a term of $b$, then other terms of $b$
should be of the form
$x_{s_{\sigma(1)}}\partial_{y_{t_1}}\wedge\cdots\wedge
x_{s_{\sigma(k)}}\partial_{y_{t_k}},(\sigma\in S_k)$.}

\vspace{0.3cm} Assume $i_{n+1}\prec_a i_{n}\prec_a\cdots\prec_a i_1$
where $\{i_1, i_2, \ldots, i_{n+1}\}$ is a permutation of
$\{0,1,\ldots,n\}$ all the time. And denote
$C_i=\{x_i\partial_{y_j}\mid 1\leq j\leq m\}, (i=0,1,\ldots,n).$

\begin{lem}
If $j<k$, then the number of elements in $C_{i_j}$ which are factors
of a term of $b$ must be equal to or grater than the number of
elements in $C_{i_k}$ which are also factors of the same term of
$b$.
\end{lem}
\begin{proof}
Suppose not. We see that one term of $b$ is of the form
$$c=x_{i_j}\partial_{y_{s_1}}\wedge
x_{i_j}\partial_{y_{s_2}}\wedge\cdots\wedge
x_{i_j}\partial_{y_{s_\alpha}}\wedge x_{i_k}\partial_{y_{t_1}}\wedge
x_{i_j}\partial_{y_{t_2}}\wedge\cdots\wedge
x_{i_j}\partial_{y_{t_\beta}}\wedge U$$, where $\alpha<\beta$ and
$U$ has no factor in $C_{i_j}$ or $C_{i_k}$. If there exist
$s_l=t_{l'}$, then we can omit the two factors
$x_{i_j}\partial_{y_{s_l}}$ and $x_{i_k}\partial_{y_{t_{l'}}}$ that
will not influence our discuss because of
$$\mbox{(if $i_j>i_k$)}\quad x_{i_j}\partial_{y_{s_l}}\wedge
x_{i_k}\partial_{y_{s_l}}\wedge\cdots\stackrel{D}{\dashrightarrow}x_{i_j}\partial_{y_{s_l}}\wedge
x_{i_k}\partial_{x_{i_j}}\wedge
x_{i_j}\partial_{y_{s_l}}\wedge\cdots(=0);$$
$$\mbox{(if $i_j<i_k$)}\quad x_{i_j}\partial_{x_{i_k}}\wedge x_{i_j}\partial_{y_{s_l}}\wedge
x_{i_k}\partial_{y_{s_l}}\wedge\cdots\stackrel{\delta}{\dashrightarrow}x_{i_j}\partial_{y_{s_l}}\wedge
x_{i_j}\partial_{y_{s_l}}\wedge\cdots(=0).$$

So we also assume that $s_l\neq t_{l'}$ for any $1\leq l\leq\alpha,
1\leq l'\leq\beta$. If $i_j>i_k$, then $x_{i_k}\partial_{x_{i_j}}$
is not a factor of $a$. Since
\begin{eqnarray*}x_{i_j}\partial_{y_{s_1}}\wedge
x_{i_j}\partial_{y_{s_2}}\wedge\cdots\wedge
x_{i_j}\partial_{y_{s_\alpha}}\wedge x_{i_k}\partial_{y_{t_1}}\wedge
x_{i_k}\partial_{y_{t_2}}\wedge\cdots\wedge
x_{i_k}\partial_{y_{t_\beta}}\wedge
\cdots\stackrel{D}{\dashrightarrow}\quad\quad\quad\\(-1)^{\alpha+l'}
x_{i_j}\partial_{y_{s_1}}\wedge\cdots\wedge
x_{i_j}\partial_{y_{s_\alpha}}\wedge
x_{i_k}\partial_{y_{t_1}}\wedge\cdots\wedge
(x_{i_k}\partial_{x_{i_j}}\wedge
x_{i_j}\partial_{y_{t_{l'}}})\wedge\cdots\wedge
x_{i_k}\partial_{y_{t_\beta}}\wedge \cdots\end{eqnarray*} and
$D(a\wedge b)=0$, there should be a term of $b$ of the form
\begin{equation}
x_{i_j}\partial_{y_{s_1}}\wedge \cdots\wedge
x_{i_j}\partial_{y_{t_{l'}}}\wedge\cdots\wedge
x_{i_j}\partial_{y_{s_\alpha}}\wedge
x_{i_k}\partial_{y_{t_1}}\wedge\cdots\wedge
x_{i_k}\partial_{y_{s_l}}\wedge\cdots\wedge
x_{i_k}\partial_{y_{t_\beta}}\wedge U
\end{equation} in which
$x_{i_j}\partial_{y_{t_{l'}}}$ and $x_{i_k}\partial_{y_{s_{l}}}$ are
at the places where $x_{i_j}\partial_{y_{s_l}}$ and
$x_{i_k}\partial_{y_{t_{l'}}}$ used to be in $c$, respectively.
However, since $\beta>\alpha$, there are not enough $l$ to match all
$l'\in\{1,2,\ldots,\beta\}$. That is impossible.

If $i_j<i_k$, then $x_{i_j}\partial_{x_{i_k}}$ is a factor of $a$.
As
\begin{eqnarray*}x_{i_j}\partial_{x_{i_k}}\wedge x_{i_j}\partial_{y_{s_1}}\wedge
x_{i_j}\partial_{y_{s_2}}\wedge\cdots\wedge
x_{i_j}\partial_{y_{s_\alpha}}\wedge x_{i_k}\partial_{y_{t_1}}\wedge
x_{i_k}\partial_{y_{t_2}}\wedge\cdots\wedge
x_{i_k}\partial_{y_{t_\beta}}\wedge
\cdots\\\stackrel{\delta}{\dashrightarrow}-
x_{i_j}\partial_{y_{s_1}}\wedge\cdots\wedge
x_{i_j}\partial_{y_{s_\alpha}}\wedge
x_{i_k}\partial_{y_{t_1}}\wedge\cdots\wedge
x_{i_j}\partial_{y_{t_{l'}}}\wedge\cdots\wedge
x_{i_k}\partial_{y_{t_\beta}}\wedge \cdots\end{eqnarray*} and
$\delta(a\wedge b)=0$, there also should be a term of $b$ of the
form
\begin{equation}
x_{i_j}\partial_{y_{s_1}}\wedge \cdots\wedge
x_{i_j}\partial_{y_{t_{l'}}}\wedge\cdots\wedge
x_{i_j}\partial_{y_{s_\alpha}}\wedge
x_{i_k}\partial_{y_{t_1}}\wedge\cdots\wedge
x_{i_k}\partial_{y_{s_l}}\wedge\cdots\wedge
x_{i_k}\partial_{y_{t_\beta}}\wedge U,
\end{equation} which is the same as we
mentioned before. There are not enough $l$ to match all
$l'\in\{1,2,\ldots,\beta\}$, either. That is impossible, too.
\end{proof}

\vspace{0.3cm} Indeed the proof of the above lemma (i.e. (4.3) and
(4.4)) also indicated other information:
\begin{lem} If $x_{i_k}\partial_{y_t}$ is a factor of a
term (denote by $c$) of $b$, then there should be other terms $c_s$
of $b$ and integers $l_s$ ($s=1,2,\dots,k-1$) such that $c_s$ comes
from $c$ by replacing $x_{i_k}\partial_{y_{t}}$ and
$x_{i_s}\partial_{y_{l_s}}$ by $x_{i_k}\partial_{y_{l_s}}$ and
$x_{i_s}\partial_{y_{t}}$, respectively. \hspace*{\fill} $\square$
\end{lem}

\vspace{0.3cm} Now we can begin our proof of the main theorem.

\noindent{\bf Proof of Theorem 4.8.} We should prove two things. One
is that the elements $a\wedge b\in\widetilde{H}(\mathcal {A}_n^m)$
can be represented as a linear combination of the elements in
$\mathcal {P}$. The other is that the elements in $\mathcal {P}$ are
linear independent.

We have discussed before that it is enough to consider the elements
satisfying condition ($\ast$). Obviously, each $a\wedge
\varphi^a_{J_1}\wedge\varphi^a_{J_2}\wedge\cdots\wedge\varphi^a_{J_t}$
satisfies ($\ast$).

Each term of $b$ can be adjusted to the ``standard'' form:
$$(x_{i_1}\partial_{y_{l_{1,1}}}\wedge x_{i_1}\partial_{y_{l_{1,2}}}
\wedge\cdots\wedge x_{i_1}\partial_{y_{l_{1,t_1}}})
\wedge\cdots\wedge (x_{i_k}\partial_{y_{l_{k,1}}}\wedge
x_{i_k}\partial_{y_{l_{k,2}}}\wedge\cdots\wedge
x_{i_k}\partial_{y_{l_{k,t_k}}}),$$ where $t_1\geq t_2\geq\cdots\geq
t_k$ (because of Lemma 4.9) and $l_{s,1}<l_{s,2}<\cdots<l_{s,t_s}
(\forall s=1,2,\ldots,k)$. For any term $c$ of $b$ of the above
``standard '' form, we define a map
$$\omega: \{\mbox{terms of}\ b\}\rightarrow \mathbb{Z}_+,\quad\quad
c\mapsto \overline{l_{1,1}l_{1,2}\cdots l_{1,t_1}\cdots
l_{k,1}l_{k,2}\cdots l_{k,t_k}},$$ where the overline means not to
multiply the elements under it but just to represent a number of
base-$(m+1)$ number system. (For example, if $l_1=1$, $l_2=12$,
$l_3=2$ and $m=19$, then $\overline{l_1l_2l_3}$ means the number
$1\times20^2+12\times20+2$.)

Since we have assumed $a\wedge b$ satisfies condition ($\ast$), each
term of $b$ should be with different values under the map $\omega$.
We call $c$ the {\em leading term} of $b$ if $\omega(c)<\omega(c')$
where $c'$ is any other term of $b$. Assume
$c=(x_{i_1}\partial_{y_{l_{1,1}}}\wedge
x_{i_1}\partial_{y_{l_{1,2}}} \wedge\cdots\wedge
x_{i_1}\partial_{y_{l_{1,t_1}}})\wedge\cdots\wedge
(x_{i_k}\partial_{y_{l_{k,1}}}\wedge
x_{i_k}\partial_{y_{l_{k,2}}}\wedge\cdots\wedge
x_{i_k}\partial_{y_{l_{k,t_k}}})$ is the leading term of $b$. We
have $l_{1,p}\leq l_{2,p}\leq\cdots\leq l_{s,p}, (\forall
s=1,2,\ldots,k; p=1,2,\ldots,t_s)$. If not, thanks to Lemma 4.10, we
may exchange certain $l_{s,p}, l_{s,p+1}$ to get another term $c'$
of $b$ such that $\omega(c')<\omega(c)$ which leads a contradiction
to our fetching way of $c$. Hence we have $t(=t_1)$ chains:
$$l_{1,p}\leq l_{2,p}\leq\cdots\leq l_{s,p}\quad (\forall s=1,2,\ldots,k;\ t_{s+1}<p\leq
t_{s}).$$

Denote $J_p=(l_{1,p},l_{2,p},\ldots,l_{s,p}),\ (\forall
s=1,2,\ldots,k;\ t_{s+1}<p\leq t_{s})$. We can observe that
$a\wedge\varphi^a_{J_1}\wedge\varphi^a_{J_2}\wedge\cdots\wedge\varphi^a_{J_t}\in\mathcal
{P}$ and $c$ is also the leading term of
$\varphi^a_{J_1}\wedge\varphi^a_{J_2}\wedge\cdots\wedge\varphi^a_{J_t}$.
Then
$a\wedge(b-\varphi^a_{J_1}\wedge\varphi^a_{J_2}\wedge\cdots\wedge\varphi^a_{J_t})$
is also an element in $\widetilde{H}(\mathcal {A}_n^m)$ and
satisfies the condition ($\ast$). Moreover, for any term $c'$ of
$(b-\varphi^a_{J_1}\wedge\varphi^a_{J_2}\wedge\cdots\wedge\varphi^a_{J_t})$,
we have $\omega(c')>\omega(c).$ Hence we can replace $b$ by
($b-\varphi^a_{J_1}\wedge\varphi^a_{J_2}\wedge\cdots\wedge\varphi^a_{J_t}$)
and use induction on $\omega(c)$. As the $\omega(c)$ becomes larger
and larger, and the $\omega(c)$ has an upper bound (because the sum
$t_1+t_2+\cdots+t_k$ is fixed), there should be an end of our
inductive process. Thus we know $a\wedge b$ can be presented as a
linear combination of elements in $\mathcal {P}$.

Now we turn to prove that $\mathcal {P}$ is a linear independent
set. Since the only change among the all terms of
$b=\varphi^a_{J_1}\wedge\varphi^a_{J_2}\wedge\cdots\wedge\varphi^a_{J_t}$
is the permutation of ${y_i}'s$, we only need to show the linear
independence of the elements in the set $$\mathcal {P}_{(t;\
k_1,k_2,\ldots,k_t)}^a=
\{\varphi^a_{J_1}\wedge\varphi^a_{J_2}\wedge\cdots\wedge\varphi^a_{J_t}\mid
a\wedge\varphi^a_{J_1}\wedge\varphi^a_{J_2}\wedge\cdots\wedge\varphi^a_{J_t}\in\mathcal
{P}, |J_s|=k_s, s=1,2,\ldots,t \},$$ where
$a\in\mathscr{B}(\mathcal{A}_{n}^0)\bigcap\widetilde{H}(\mathcal
{A}_{n}^0)$ and $t, k_1,k_2,\ldots,k_t\in\mathbb{Z}_+$ with $k_1\geq
k_2\geq\cdots\geq k_t$ are fixed. ($|J_s|$ means the length of
$J_s$. For example, if $J_s=(j_1,j_2,\ldots,j_k)$, then $|J_s|=k$.)

For any two elements of ${P}_t^a$, their leading terms are different
from each other. So are the values of these two leading terms under
the map $\omega$.

If $\sum_{s=1}^p\alpha_sb_s=0$ where $0\neq\alpha_s\in\mathbb{F}$
and $0\neq b_s\in\mathcal {P}_t^a$ with leading term $c_s$. Suppose
$\omega(c_{s_1})<\omega(c_{s_2})<\cdots<\omega(c_{s_p})$, then the
leading term of $\sum_{s=1}^p\alpha_sb_s$ is $c_{s_1}\neq 0$. But
the leading term of $0$ is of course $0$. Contradiction! Therefore
the elements in $\mathcal {P}_{(t;\ k_1,k_2,\ldots,k_t)}^a$ are
linear independent. So are the elements in $\mathcal {P}$.
\hspace*{\fill} $\square$

\vspace{0.3cm} By Theorem 4.8, the calculation of the Betti numbers
can be transformed to a combinatorial problem. Precisely, we should
to computer the number of the elements in a such set:
\begin{eqnarray*}
&&\mathcal
{S}_{m,n}^i=\{[(l_{1,1},l_{1,2},\ldots,l_{1,t_1}),(l_{2,1},l_{2,2},\ldots,l_{2,t_2}),
\ldots,(l_{k,1},l_{k,2},\ldots,l_{k,t_k})]\mid\\&&\quad\quad\quad\quad1\leq
k\leq n+1,\ t_k\leq t_{k-1}\leq\cdots\leq t_1,\
t_1+t_2+\cdots+t_k=i,\\&&\quad\quad\quad\quad1\leq
l_{s,1}<l_{s,2}<\cdots<l_{s,t_s}\leq m,\ l_{1,p}\leq
l_{2,p}\leq\cdots\leq l_{s,p},1\leq s\leq k, p\leq t_s\},
\end{eqnarray*}
in which $0\leq i\leq mn+m$.

Using the definitions and notations of Young tableaux in [Fw], we
say $\mathcal {S}_{m,n}^i$ is the set of all Young tableaux whose
entries are taken from $\{1,2,\ldots,m\}$ and whose shape is
$\lambda=(\lambda_1\geq\lambda_2\geq\cdots\geq\lambda_s)\vdash i$
with $\lambda_1\leq n+1$ and $s\in\mathbb{Z}_+$.

Denote by $d_\lambda(m)$ the number of Young tableaux on the shape
$\lambda$ with entries in $\{1,2,\ldots,m\}$. One has $$|\mathcal
{S}_{m,n}^i|=\sum_{\lambda=(\lambda_1\geq\lambda_2\geq\cdots\geq\lambda_s)\vdash
i;\ \lambda_1\leq n+1}d_\lambda(m).$$

There is a hook length formula for the number $d_\lambda(m)$ due to
Stanley (c.f. [Fw]):
$$d_\lambda(m)=\prod_{(i,j)\in\lambda}\frac{m+j-i}{h_\lambda(i,j)},$$
where $h_\lambda(i,j)$ is the hook length in the $i$-th row and
$j$-th column of shape $\lambda$.

The following corollary can be obtained by Theorem 4.8 and Corollary
4.5.
\begin{cor}
The generating function of the Betti numbers of $\mathcal {A}_n^m$
is $$\sum_{i=0}^\infty b_it^i=(\sum_{i=0}^{mn+m}|\mathcal
{S}_{m,n}^i|t^i)\frac{\prod_{j=1}^n(1-t^{j+1})}{(1-t)^n},$$ where
$|\mathcal
{S}_{m,n}^i|=\sum\limits_{\lambda=(\lambda_1\geq\lambda_2\geq\cdots\geq\lambda_s)\vdash
i;\ \lambda_1\leq
n+1}\prod_{(i,j)\in\lambda}\frac{m+j-i}{h_\lambda(i,j)}$.
\hspace*{\fill} $\square$
\end{cor}

\vspace{0.5cm} By the total order (4.1) and Theorem 4.4, there is a
one-to-one correspondence between
$\mathscr{B}(\mathcal{A}_{n}^0)\bigcap \widetilde{H}(\mathcal
{A}_{n}^0)$ and the $(n+1)$-th symmetric group $S_{n+1}$. Precisely,
the element which corresponds to $\sigma\in S_{n+1}$ belongs to
$H^k_\tau(\mathcal {A}_n^0)$ with $(-1)^k=\mbox{sign}(\sigma)$ and
$\tau=\sum\limits_{0\leq i<j\leq n;\
\sigma(i)>\sigma(j)}\epsilon_{\sigma(i)}-\epsilon_{\sigma(j)}$.

Apply Euler-Poincar$\acute{\mbox{e}}$ Principle to $\mathcal
{A}_n^0$. We get an identity:
\begin{equation}
\prod_{0\leq i<j\leq n}(1-e^{\epsilon_j-\epsilon_i})=\sum_{\sigma\in
S_{n+1}}\mbox{sign}(\sigma)\prod_{0\leq i<j\leq n;
\sigma(i)>\sigma(j)}e^{\epsilon_{\sigma(i)}-\epsilon_{\sigma(j)}}.
\end{equation}
Multiply the both sides of (4.5) by
$\prod_{i=0}^ne^{(n-i)\epsilon_i}$. We get
\begin{equation*}
\prod_{0\leq i<j\leq
n}(e^{\epsilon_i}-e^{\epsilon_j})=\sum_{\sigma\in
S_{n+1}}\mbox{sign}(\sigma)\prod_{i=0}^n
e^{(n-i)\epsilon_{\sigma(i)}},
\end{equation*}
which is just the Vandermonde Determinant!

For $\mathcal {A}_n^m$, there is a one-to-one correspondence between
$\mathcal {P}$ and $\{(\sigma, T)\mid \sigma\in S_{n+1},
T\in{S}_{m,n}^t, 0\leq t\leq mn+m\}$. Now the
Euler-Poincar$\acute{\mbox{e}}$ Principle induces an identity as
follows:
\begin{eqnarray*}
&&\prod_{0\leq i\leq n;\ i<j\leq m+n}(1-e^{\epsilon_j-\epsilon_i})=
\\&&\sum_{\sigma\in
S_{n+1}}\mbox{sign}(\sigma)\prod_{0\leq i<j\leq n;
\sigma(i)>\sigma(j)}e^{\epsilon_{\sigma(i)}-\epsilon_{\sigma(j)}}
(\sum_{t=0}^{mn+m}\sum_{T\in{S}_{m,n}^t}
(-1)^t(\prod_{k=1}^me^{c_k\epsilon_{n+k}})(\prod_{k=1}^s
e^{-\lambda_k\epsilon_{\sigma(n-k+1)}})),
\end{eqnarray*}
where $(\lambda_1\geq\lambda_2\geq\cdots\geq\lambda_s)\vdash t$ and
$(c_1,c_2,\cdots,c_m)$ are the shape and content of $T$,
respectively.

\section{Final Remarks about $H(L_{1}(\mathcal
{T}^{d}))$}

In [L], we also introduced a class of solvable Lie algebras
$L_{1}(\mathcal {T}^{d})$, which is an extension of $L_{0}(\mathcal
{T}^{d})$ with $\mathcal {H}=\{x_i\partial_{x_i}\mid \iota_i\in
\mathcal {N}\}$. The cohomology of $L_{1}(\mathcal {T}^{d})$ is much
simpler than that of $L_{0}(\mathcal {T}^{d})$. In fact, using the
following lemma, we can obtain a theorem about $H(L_{1}(\mathcal
{T}^{d}))$ immediately.

To describe the Lemma, we have to introduce some notations and
definitions firstly.

Given a finite dimensional Lie algebra $\mathcal {G}$, suppose
$g_1,g_2,...,g_t\in\mathcal {G}$ are pairwise commuting elements
such that $\mathcal {G}$ possesses a basis consisting of the vectors
which are eigenvectors for all the operators $ad\
g_i:g\mapsto[g_i,g]$. Denote $$\mathcal {G}_{(\lambda_1,
\lambda_2,\ldots, \lambda_t)}=\{g\in\mathcal
{G}\mid[g_i,g]=\lambda_ig,i=1,2,\ldots, t\}.$$ It is obvious that
$$[\mathcal {G}_{(\lambda_1,\ldots,\lambda_t)},\mathcal {G}_{(\mu_1,\dots,\mu_t)}]
\subset\mathcal {G}_{(\lambda_1+\mu_1,\ldots,\lambda_t+\mu_t)}.$$
Denote $$\wedge^k_{(\lambda_1,\lambda_2,\ldots,\lambda_t)}\mathcal
{G}=\mbox{Span}\{r_1\wedge r_2\wedge\cdots\wedge r_k\mid
r_i\in\mathcal {G}_{(\lambda_{i_1},\ldots,\lambda_{i_t})}(1\leq
i\leq k),\sum_{i=1}^k\lambda_{i_j}=\lambda_j(1\leq j\leq t)\}$$ and
$$\wedge_{(\lambda_1,\lambda_2,\ldots,\lambda_t)}\mathcal
{G}=\bigoplus_{k\geq0}\wedge^k_{(\lambda_1,\lambda_2,\ldots,\lambda_t)}\mathcal
{G}.$$
\begin{lem}
{\em\bf ([F])} The inclusion $\wedge_{(0, 0, \ldots, 0)}\mathcal
{G}\rightarrow\wedge \mathcal {G}$ induces an isomorphism in
cohomology.\hfill$\Box$
\end{lem}

Hence there comes a theorem:

\begin{thm}
The cohomology group of $L_{1}(\mathcal {T}^{d})$ is isomorphic to
$\wedge\mathcal {H}=\oplus_{i\geq0}\wedge^i\mathcal {H}$, where
$\mathcal {H}=\{x_i\partial_{x_i}\mid \iota_i\in \mathcal {N}\}$.
\end{thm}
\begin{proof}
Take the natural basis of $L_{0}(\mathcal {T}^{d})$, i.e.
$$B(\mathcal
{T}^d)=\{\partial_{x_i},(\prod_{\iota_s\in \mathcal
{C}_j\backslash\{\iota_j\}}x_s^{m_s})\partial_{x_j}\mid
\iota_i\in\Lambda, m_s\in\mathbb{N},\sum_{\iota_s\in \mathcal
{C}_j\backslash\{\iota_j\}}m_s\kappa_{s,j}\leq \kappa_{j}\}.$$ Then
$B(\mathcal {T}^d)\bigcup\mathcal {H}$ is a basis of $L_{1}(\mathcal
{T}^{d})$. Furthermore, the elements in $\mathcal {H}$ are pairwise
commutative and the elements in $B(\mathcal {T}^d)\bigcup\mathcal
{H}$ are the eigenvectors for all the operators
$ad(x_i\partial_{x_i}),(\forall x_i\partial_{x_i}\in\mathcal {H})$.
Hence we can define the $\mathcal {G}_{(\lambda_1, \lambda_2,\ldots,
\lambda_t)}$ and
$\wedge_{(\lambda_1,\lambda_2,\ldots,\lambda_t)}\mathcal {G}$ for
$\mathcal {G}=L_{1}(\mathcal {T}^{d})$.

Take any element $r_1\wedge r_2\wedge\cdots\wedge r_p\in \wedge_{(0,
0, \ldots, 0)}\mathcal {G},(r_i\in B(\mathcal {T}^d)\bigcup\mathcal
{H})$. If $r_{i_1}=f_1\partial_{x_{s_1}}\in B(\mathcal {T}^d)$, then
we have $[x_{s_1}\partial_{x_{s_1}},r_{i_1}]=-r_{i_1}$. So there
must be an $r_{i_2}$ and $\alpha>0$ such that
$[x_{s_1}\partial_{x_{s_1}},r_{i_2}]=\alpha r_{i_2}$. Thus we can
assume $r_{i_2}=f_2\partial_{x_{s_2}}\in B(\mathcal {T}^d)$
satisfying that $x_{s_1}$ is a factor of $f_2$ and $\iota_{s_2}\in
\mathcal {D}_{s_1}\backslash\{\iota_{s_1}\}$. Analogously, there
must be $r_{i_k}=f_k\partial_{x_{s_k}}$ such that $x_{s_{k-1}}$ is a
factor of $f_k$ and $\iota_{s_k}\in \mathcal
{D}_{s_{k-1}}\backslash\{\iota_{s_{k-1}}\}$. $k$ can be taken as
arbitrary positive integer, but there exists no infinite chain
$\iota_{s_1},\iota_{s_2},\ldots,\iota_{s_i},\ldots$ such that
$\iota_{s_i}\in \mathcal
{D}_{s_{i-1}}\backslash\{\iota_{s_{i-1}}\}$. Thus
$r_1,r_2,\ldots,r_p$ must be all in $\mathcal {H}$, i.e.
$\wedge_{(0, 0, \ldots, 0)}\mathcal {G}\subset\wedge\mathcal {H}$.
As it is obvious that $\wedge\mathcal {H}\subset\wedge_{(0, 0,
\ldots, 0)}\mathcal {G}$, we get $\wedge\mathcal {H}=\wedge_{(0, 0,
\ldots, 0)}\mathcal {G}$.

Thanks to Lemma 5.1, $\wedge\mathcal {H}$ is isomorphic to the
cohomology group of $L_{1}(\mathcal {T}^{d})$.
\end{proof}

\vspace{1cm} \noindent{\Large \bf Acknowledgements}

\vspace{0.3cm} Thanks are due to my thesis advisor Prof. Xiaoping Xu
for his help not only in my study, but also in my daily life.

I also thank the referee for helpful suggestions which improved the
exposition of this paper.

\vspace{1cm}

\noindent{\Large \bf References}

\hspace{0.5cm}

\begin{description}

\item[{[ACJ]}] G.F.Armstrong, G.Cairns and B.Jessup, Explicit Betti
numbers for a family of nilpotent Lie algebras, {\it Proc. Amer.
Math. Soc.} {\bf 125} (1997), 381-385.

\item[{[AS]}] G.F.Armstrong and S.Sigg, On the cohomology of a class of
nilpotent Lie algebras, {\it Bull. Austral. Math. Soc.} {\bf 54(2)}
(1996), 517-527.

\item[{[B]}] M.Bordemann, Nondegenerate invariant bilinear forms on
nonassociative algebras, {\it Acta. Math. Univ. Comenian.} {\bf
66(2)} (1997), 151-201.

\item[{[Br]}] R.Bott, Homogeneous vector bundles, {\it Ann. of
Math.} {\bf 66} (1957), 203-248.

\item[{[BGG]}] I.N.Bernstein, I.M.Gelfand and S.I.Gelfand,
Differential operators on the base affine space and a study of
$\mathcal {G}$-modules, {\it Lie groups and their representations.}
Summer School in Group Representations. Bolyai Janos Math. Soc.,
Budapest 1971, pp.21-64, New York: Halsted 1975.

\item[{[CJP]}] G.Cairns, B.Jessup and J.Pitkethly, On the Betti
numbers of nilpotent Lie algebras of small dimension, {\it
Integrable Systems and Foliations}, Birkh\"{a}user, Boston, MA,
1997, 19-31.

\item[{[D]}] J.Dixmier, Cohomologies des alg$\grave{\mbox{e}}$bres de Lie nilpotents,
{\it Acta Sci. Math. Szeged} {\bf 16} (1955), 246-250.

\item[{[DS]}] Ch.Deninger and W.Singhof, On the cohomology of nilpotent
Lie algebras, {\it Bull. Soc. Math. France} {\bf 116} (1988), 3-14.

\item[{[F]}] D.B.Fuchs, {\it Cohomology of the Infinite Dimensional Lie
Algebras}, Consultant Bureau, New York, 1987.

\item[{[FM]}] A.Fialowski and D.Millionschikov, Cohomology of graded
Lie algebras of maximal class, {\it J.Algebra} {\bf 296}
(2006),157-176.

\item[{[Fw]}] W.Fulton, {\it Young Tableaux: with Applications to Representation Theory and Geometry},
LMSST 35, Cambridge University Press, 1997.

\item[{[GL]}] H.Garland and J.Lepowsky, Lie algebra homology and the
Macdonald-Kac formulas, {\it Inv. Math.} {\bf 34} (1976), 37-76.

\item[{[H]}] J.E.Humphreys, {\it Reflection Groups and Coxeter
Groups}, Cambridge University Press, 1990.

\item[{[K]}] B.Kostant, Lie algebra cohomology and the generalized
Borel-Weil theorem, {\it Ann. of Math.} {\bf 74} (1961), 329-387.

\item[{[KK]}] S.J.Kang and M.H.Kim, {\it Dimension Formula for Graded Lie Algebras and its
Applications}, RIM-GARC Preprint Series 96\_\ 51, 1996.

\item[{[L]}] L.Luo, Oriented tree diagram Lie algebras and their abelian
ideals, submitted.

\item[{[S]}] L.J.Santharoubane, Cohomology of Heisenberg Lie
algebras, {\it Proc. Amer. Math. Soc.} {\bf 87} (1983), 23-28.

\item[{[X]}] X.Xu, Tree diagram Lie algebras of differential
operators and evolution partial differential equations, {\it J. Lie
Theory} {\bf 16} (2006), 691-718.

\end{description}

\end{document}